\documentclass[12pt]{article}
\usepackage{latexsym,amsfonts,amssymb}
\usepackage[dvips]{color}
\usepackage{colordvi,multicol}
\usepackage{amsmath}
\usepackage{graphicx}
\usepackage{pdflscape}
\usepackage{rotating}

\usepackage{latexsym,amsfonts,amssymb}
\usepackage[dvips]{color}
\usepackage{colordvi,multicol}
\usepackage{amsmath}
\usepackage{graphicx}
\usepackage{pdflscape}
\usepackage{rotating}
\usepackage{breqn}

\def \cal{\mathcal}

\newtheorem{thm}{Theorem}[section]
\newtheorem{cor}[thm]{Corollary}
\newtheorem{lem}[thm]{Lemma}

\newtheorem{rem}[thm]{Remark}

\date{}
\begin{document}

\title{\bf Modulus representation of the Riemann
$\xi$ function and polynomial inequalities equivalent to the Riemann hypothesis}
\author{Wei Sun\\ \\ \\
  {\small Department of Mathematics and Statistics}\\
{\small Concordia University}\\
{\small Montreal, H3G 1M8,  Canada}\\
{\small  wei.sun@concordia.ca}}

\maketitle

\begin{abstract}
\noindent We use the Jacobi theta function to give a representation of the modulus of the Riemann $\xi$ function. Based on this modulus representation,  we show that the Riemann hypothesis is equivalent to the validity of a family of polynomial inequalities. We also present some preliminary results on the polynomial inequalities.
\end{abstract}



\section{Introduction and the main results}
The Riemann hypothesis (RH) is the conjecture that the Riemann zeta function $\zeta(s)$ has its nontrivial zeros only on the critical line ${\rm Re}(s)=\frac{1}{2}$. Let $\xi(s)$ be the Riemann xi function, i.e.,
$$
\xi(s)=\frac{1}{2}s(s-1)\pi^{-\frac{s}{2}}\Gamma\left(\frac{s}{2}\right)\zeta(s),\ \ \ \ s\in\mathbb{C},
$$
where $\Gamma(s)$ is the Gamma function. Then, the RH is equivalent to the statement that the zeros of $\xi(s)$  are all located on the critical line. We refer the reader to Edwards's book \cite{E} for the terminology used in this paper, and to Wikipedia \cite{W} and references therein for an introduction to the RH.

Let $\psi$ be the Jacobi theta function:
$$
\psi(y)=\sum_{n=1}^{\infty}e^{-\pi n^2y},\ \ \ \ y>0.
$$
It is well known that  $\xi$ can be represented as
\begin{eqnarray}\label{Rie}
\xi(s)=\frac{1}{2}+\frac{s(s-1)}{2}\int_1^{\infty}\left(y^{\frac{s}{2}-1}+y^{-\frac{s+1}{2}}\right)\psi(y)dy,\ \ \ \ s\in\mathbb{C}.
\end{eqnarray}
See Riemann's memoir \cite{R} and \cite[Pages 16 and 17]{E}.

In this paper, we will first prove the following representation of the modulus of $\xi$.
\begin{thm}\label{thm1}
For any $\tau\in\mathbb{R}$ and $t\in\mathbb{R}$,
\begin{eqnarray}\label{Apr4v}
&&4\left|\xi\left(\frac{1}{2}+\tau-it\right)\right|^2\nonumber\\
&=&2\left[\left(t^2+\tau^2+\frac{1}{4}\right)^2-\tau^2\right]\int_{0}^{\infty}\int_{0}^{\infty}\cos\left(\frac{t}{2}\ln (uv)\right)u^{\frac{2\tau-3}{4}}v^{\frac{-2\tau-3}{4}}1_{\{uv>1\}}\psi(u)\psi(v)dudv\nonumber\\
&&-t^2\int_1^{\infty}\left[(1+2\tau)(u^{\tau-\frac{1}{2}}+u^{-\tau-1})+(1-2\tau)(u^{-\tau-\frac{1}{2}}+u^{\tau-1})\right]\psi(u)du\nonumber\\
&&-2\left(\tau^2-\frac{1}{4}\right)^2\int_{0}^{\infty}\int_{0}^{\infty}u^{\frac{2\tau-3}{4}}v^{\frac{-2\tau-3}{4}}1_{\{uv>1\}}\psi(u)\psi(v)dudv\nonumber\\
&&+\left[1+\left(\tau^2-\frac{1}{4}\right)\int_1^{\infty}\left(u^{\frac{2\tau-3}{4}}+u^{\frac{-2\tau-3}{4}}\right)\psi(u)du\right]^2.
\end{eqnarray}
\end{thm}

According to Sondow and Dumitrescu \cite[Corollary 1]{SD}, the RH is equivalent to the statement that for any fixed $t\in\mathbb{R}$ the function $\tau\mapsto|\xi(\frac{1}{2}+\tau-it)|$ is increasing for $\tau>0$. Based on this monotonicity property and (\ref{Apr4v}), we are able to show that the RH is equivalent to the validity of a family of polynomial inequalities.

For $\tau\in\mathbb{R}$, define

\noindent (i) \begin{eqnarray}\label{July4b}
a_{\tau}(0)&:=&\left[1-\left(\frac{1}{4}-\tau^2\right)\int_1^{\infty}\left(u^{\frac{2\tau-3}{4}}+u^{\frac{-2\tau-3}{4}}\right)\psi(u)du\right]\nonumber\\
&&\cdot\left[4\tau\int_1^{\infty}\left(u^{\frac{2\tau-3}{4}}+u^{\frac{-2\tau-3}{4}}\right)\psi(u)du-\left(\frac{1}{4}-\tau^2\right)\int_1^{\infty}\left(u^{\frac{2\tau-3}{4}}-u^{\frac{-2\tau-3}{4}}\right)(\ln u)\psi(u)du\right],\nonumber\\
&&
\end{eqnarray}
\begin{eqnarray*}
a_{\tau}(1)&:=&8\tau\int_{0}^{\infty}\int_{0}^{\infty}u^{\frac{2\tau-3}{4}}v^{\frac{-2\tau-3}{4}}1_{\{uv>1\}}\psi(u)\psi(v)dudv\\
&&+\tau\left(\frac{1}{4}-\tau^2\right)\int_{0}^{\infty}\int_{0}^{\infty}\left[\ln (uv)\right]^{2}u^{\frac{2\tau-3}{4}}v^{\frac{-2\tau-3}{4}}1_{\{uv>1\}}\psi(u)\psi(v)dudv\\
&&+2\left(\frac{1}{4}+\tau^2\right)\int_{0}^{\infty}\int_{0}^{\infty}\ln\left(\frac{u}{v}\right)u^{\frac{2\tau-3}{4}}v^{\frac{-2\tau-3}{4}}1_{\{uv>1\}}\psi(u)\psi(v)dudv\\
&&-\frac{(\frac{1}{4}-\tau^2)^2}{8}\int_{0}^{\infty}\int_{0}^{\infty}\left[\ln (uv)\right]^{2}\ln\left(\frac{u}{v}\right)u^{\frac{2\tau-3}{4}}v^{\frac{-2\tau-3}{4}}1_{\{uv>1\}}\psi(u)\psi(v)dudv\\
&&-2\int_1^{\infty}\left[(u^{\tau-\frac{1}{2}}+u^{-\tau-1})-(u^{-\tau-\frac{1}{2}}+u^{\tau-1})\right]\psi(u)du\nonumber\\
&&-\int_1^{\infty}\left[(1+2\tau)(u^{\tau-\frac{1}{2}}-u^{-\tau-1})+(1-2\tau)(u^{\tau-1}-u^{-\tau-\frac{1}{2}})\right](\ln u)\psi(u)du,
\end{eqnarray*}
and for $k\ge2$,
\begin{eqnarray*}
a_{\tau}(k)
&:=&-\frac{(-1)^k\tau}{(2k-2)!\cdot 2^{2k-5}}\int_{0}^{\infty}\int_{0}^{\infty}\left[\ln (uv)\right]^{2k-2}u^{\frac{2\tau-3}{4}}v^{\frac{-2\tau-3}{4}}1_{\{uv>1\}}\psi(u)\psi(v)dudv\\
&&-\frac{(-1)^k\tau(\frac{1}{4}-\tau^2)}{(2k)!\cdot 2^{2k-3}}\int_{0}^{\infty}\int_{0}^{\infty}\left[\ln (uv)\right]^{2k}u^{\frac{2\tau-3}{4}}v^{\frac{-2\tau-3}{4}}1_{\{uv>1\}}\psi(u)\psi(v)dudv\\
&&+\frac{(-1)^k}{(2k-4)!\cdot 2^{2k-4}}\int_{0}^{\infty}\int_{0}^{\infty}\left[\ln (uv)\right]^{2k-4}\ln\left(\frac{u}{v}\right)u^{\frac{2\tau-3}{4}}v^{\frac{-2\tau-3}{4}}1_{\{uv>1\}}\psi(u)\psi(v)dudv\\
&&-\frac{(-1)^k(\frac{1}{4}+\tau^2)}{(2k-2)!\cdot 2^{2k-3}}\int_{0}^{\infty}\int_{0}^{\infty}\left[\ln (uv)\right]^{2k-2}\ln\left(\frac{u}{v}\right)u^{\frac{2\tau-3}{4}}v^{\frac{-2\tau-3}{4}}1_{\{uv>1\}}\psi(u)\psi(v)dudv\\
&&+\frac{(-1)^k(\frac{1}{4}-\tau^2)^2}{(2k)!\cdot 2^{2k}}\int_{0}^{\infty}\int_{0}^{\infty}\left[\ln (uv)\right]^{2k}\ln\left(\frac{u}{v}\right)u^{\frac{2\tau-3}{4}}v^{\frac{-2\tau-3}{4}}1_{\{uv>1\}}\psi(u)\psi(v)dudv.
\end{eqnarray*}

\noindent (ii) 
\begin{eqnarray*}
a_{\tau,1}(1)&:=&8\tau\int_{0}^{\infty}\int_{0}^{\infty}u^{\frac{2\tau-3}{4}}v^{\frac{-2\tau-3}{4}}1_{\{uv>1\}}\psi(u)\psi(v)dudv\\
&&+\tau\left(\frac{1}{4}-\tau^2\right)\int_{0}^{\infty}\int_{0}^{\infty}\left[\ln (uv)\right]^{2}u^{\frac{2\tau-3}{4}}v^{\frac{-2\tau-3}{4}}1_{\{uv>1\}}\psi(u)\psi(v)dudv\\
&&+2\left(\frac{1}{4}+\tau^2\right)\int_{0}^{\infty}\int_{0}^{\infty}\ln\left(\frac{u}{v}\right)u^{\frac{2\tau-3}{4}}v^{\frac{-2\tau-3}{4}}1_{\{uv>1\}}\psi(u)\psi(v)dudv\\
&&-2\int_1^{\infty}\left[(u^{\tau-\frac{1}{2}}+u^{-\tau-1})-(u^{-\tau-\frac{1}{2}}+u^{\tau-1})\right]\psi(u)du\nonumber\\
&&-\int_1^{\infty}\left[(1+2\tau)(u^{\tau-\frac{1}{2}}-u^{-\tau-1})+(1-2\tau)(u^{\tau-1}-u^{-\tau-\frac{1}{2}})\right](\ln u)\psi(u)du,
\end{eqnarray*}
for $n\ge 2$,
\begin{eqnarray*}
a_{\tau,n}(2n-1)
&:=&\frac{\tau}{(4n-4)!\cdot 2^{4n-7}}\int_{0}^{\infty}\int_{0}^{\infty}\left[\ln (uv)\right]^{4n-4}u^{\frac{2\tau-3}{4}}v^{\frac{-2\tau-3}{4}}1_{\{uv>1\}}\psi(u)\psi(v)dudv\\
&&+\frac{\tau(\frac{1}{4}-\tau^2)}{(4n-2)!\cdot 2^{4n-5}}\int_{0}^{\infty}\int_{0}^{\infty}\left[\ln (uv)\right]^{4n-2}u^{\frac{2\tau-3}{4}}v^{\frac{-2\tau-3}{4}}1_{\{uv>1\}}\psi(u)\psi(v)dudv\\
&&-\frac{1}{(4n-6)!\cdot 2^{4n-6}}\int_{0}^{\infty}\int_{0}^{\infty}\left[\ln (uv)\right]^{4n-6}\ln\left(\frac{u}{v}\right)u^{\frac{2\tau-3}{4}}v^{\frac{-2\tau-3}{4}}1_{\{uv>1\}}\psi(u)\psi(v)dudv\\
&&+\frac{\frac{1}{4}+\tau^2}{(4n-4)!\cdot 2^{4n-5}}\int_{0}^{\infty}\int_{0}^{\infty}\left[\ln (uv)\right]^{4n-4}\ln\left(\frac{u}{v}\right)u^{\frac{2\tau-3}{4}}v^{\frac{-2\tau-3}{4}}1_{\{uv>1\}}\psi(u)\psi(v)dudv,
\end{eqnarray*}
and for $n\in\mathbb{N}$,
\begin{eqnarray*}
a_{\tau,n}(2n):=\frac{1}{(4n-4)!\cdot 2^{4n-4}}\int_{0}^{\infty}\int_{0}^{\infty}\left[\ln (uv)\right]^{4n-4}\ln\left(\frac{u}{v}\right)u^{\frac{2\tau-3}{4}}v^{\frac{-2\tau-3}{4}}1_{\{uv>1\}}\psi(u)\psi(v)dudv.
\end{eqnarray*}

\noindent (iii) for $n\in\mathbb{N}$,
$$
f_{\tau,n}(t):=\left[\sum_{k=0}^{2n-2}a_{\tau}(k)t^{2k}\right]+a_{\tau,n}(2n-1)t^{4n-2}+a_{\tau,n}(2n)t^{4n},\ \ \ \ t\in\mathbb{R}.
$$

\begin{thm}\label{thm2} (i) If the RH is true, then $f_{\tau,n}(t)>0$ for any $t\in\mathbb{R}$, $\tau\in(0,\infty)$ and $n\in\mathbb{N}$.

\noindent (ii) If there exists an increasing sequence of natural numbers $\{n_k\}$ such that $f_{\tau,n_k}(t)\ge0$ for any $t\in\mathbb{R}$, $\tau\in(0,\frac{1}{2})$ and $k\in\mathbb{N}$, then the RH is true.
\end{thm}

Theorem \ref{thm2} transforms the validity of the RH into the problem of checking nonnegativity of univariate polynomials, which can be explored using various approaches including e.g. semidefinite programming. There are different methods to determine the number of real roots of a general univariate polynomial with real-coefficients. We refer the reader to Yang et al. \cite{Y1} (cf. also \cite{Y2}) for an explicit algorithm using the Sylvester matrix and to \cite{P} for using the signature of the Hermite form.

For a polynomial $f$ on $\mathbb{R}$, let Discr$(f)$ denote its discriminant and ${\cal N}(f)$ denote its number of real roots.  As a consequence of Theorem \ref{thm2}, we can obtain the following necessary and sufficient condition for the validity of the RH.

\begin{thm}\label{pro99}
The RH is true if and only if there exists an increasing sequence of natural numbers $\{n_k\}$ such that for any $k\in\mathbb{N}$ and $\tau\in(0,\frac{1}{2})$,
\begin{equation}\label{Nov1a}
{\rm Discr}(f_{\tau,n_k})=0\Longrightarrow {\cal N}(f_{\tau,n_k})=0.
\end{equation}
\end{thm}

In addition, we have the following result.

\begin{cor}\label{cor111}
(i) If the RH is true, then Discr$(f_{\tau,n})\ge 0$ for any $\tau\in(0,\infty)$ and $n\in\mathbb{N}$.

\noindent (ii) If there exists an increasing sequence of natural numbers $\{n_k\}$ such that  
\begin{equation}\label{Nov4a}
\lim_{k\rightarrow\infty}\sup\left\{\tau\in\left[0,\frac{1}{2}\right):{\rm Discr}(f_{\tau,n_k})=0\right\}=0,
\end{equation}
then the RH is true.

\noindent (iii) For any $\tau\in(0,\frac{1}{2})$, Discr$(f_{\tau,1})\ge 0$, and Discr$(f_{\tau,1})= 0\Longleftrightarrow |a_{\tau,1}(1)|=2\sqrt{a_{\tau,1}(2)a_{\tau}(0)}$.

\noindent (iv) The polynomial $f_{\tau,1}(t)>0$ for any $t\in\mathbb{R}$ and $\tau\in(0,\frac{1}{2})$ if and only if  $a_{\tau,1}(1)\not=-2\sqrt{a_{\tau,1}(2)a_{\tau}(0)}$ for any $\tau\in(0,\frac{1}{2})$. Hence, the RH is not true if there exists $\tau\in(0,\frac{1}{2})$ such that $a_{\tau,1}(1)=-2\sqrt{a_{\tau,1}(2)a_{\tau}(0)}$.
\end{cor}

The rest of the paper is organized as follows. We give the proofs of Theorems \ref{thm1} and \ref{thm2} in Sections 2 and 3, respectively. The proofs of  Theorem \ref{pro99} and Corollary \ref{cor111} will be given in Section 4. 

\section{Proof of Theorem \ref{thm1}}\setcounter{equation}{0}

We assume without loss of generality that $\tau,t>0$. Define
$$
{\cal F}_{\tau,t}(z):=\frac{\xi(\tau+\frac{1}{2}+z)}{(\tau+\frac{1}{2}+z)(\tau-\frac{1}{2}+z)}\cdot \frac{\xi(-\tau+\frac{1}{2}+z)}{(-\tau+\frac{1}{2}+z)(-\tau-\frac{1}{2}+z)},\ \ \ \ z\in\mathbb{C}.
$$
Let $c>\frac{2\tau+3}{4}$. By Cauchy's integral formula, we get
\begin{eqnarray}\label{July5a}
&&\frac{1}{2\pi i }\int_{c-i\infty}^{c+i\infty}\frac{z{\cal F}_{\tau,t}(z)}{(z+it)(z-it)}dz+\frac{1}{2\pi i }\int_{-c+i\infty}^{-c-i\infty}\frac{z{\cal F}_{\tau,t}(z)}{(z+it)(z-it)}dz\nonumber\\
&=&{\cal F}_{\tau,t}(it)+\frac{1}{4\tau}\Bigg\{\left[\frac{\xi(2\tau+1)\xi(1)}{t^2+(\tau+\frac{1}{2})^2}+\frac{\xi(0)\xi(-2\tau)}{t^2+(\tau+\frac{1}{2})^2}\right]-\left[\frac{\xi(2\tau)\xi(0)}{t^2+(\tau-\frac{1}{2})^2}+\frac{\xi(1)\xi(-2\tau+1)}{t^2+(\tau-\frac{1}{2})^2}\right]\Bigg\}\nonumber\\
&=&{\cal F}_{\tau,t}(it)+\frac{1}{2\tau }\Bigg\{\frac{\xi(0)\xi(-2\tau)}{t^2+(\tau+\frac{1}{2})^2}-\frac{\xi(0)\xi(2\tau)}{t^2+(\tau-\frac{1}{2})^2}\Bigg\}\nonumber\\
&=&{\cal F}_{\tau,t}(it)+\frac{[\xi(-2\tau)-\xi(2\tau)]t^2+[(\tau-\frac{1}{2})^2\xi(-2\tau)-(\tau+\frac{1}{2})^2\xi(2\tau)]}{4\tau[t^2+(\tau+\frac{1}{2})^2][t^2+(\tau-\frac{1}{2})^2]}.
\end{eqnarray}
Note that
\begin{eqnarray*}
\frac{1}{\pi i }\int_{c-i\infty}^{c+i\infty}\frac{z{\cal F}_{\tau,t}(z)}{(z+it)(z-it)}dz
&=&\frac{1}{\pi i }\int_{-c+i\infty}^{-c-i\infty}\frac{y{\cal F}_{\tau,t}(-y)}{(-y+it)(-y-it)}dy\\
&=&\frac{1}{\pi i }\int_{-c+i\infty}^{-c-i\infty}\frac{z{\cal F}_{\tau,t}(z)}{(z+it)(z-it)}dz.
\end{eqnarray*}
Then, by (\ref{July5a}), we obtain that
\begin{eqnarray}\label{July3d}
&&{\cal F}_{\tau,t}(it)+\frac{[\xi(-2\tau)-\xi(2\tau)]t^2+[(\tau-\frac{1}{2})^2\xi(-2\tau)-(\tau+\frac{1}{2})^2\xi(2\tau)]}{4\tau[t^2+(\tau+\frac{1}{2})^2][t^2+(\tau-\frac{1}{2})^2]}\nonumber\\
&=&\frac{1}{\pi i }\int_{c-i\infty}^{c+i\infty}\frac{z{\cal F}_{\tau,t}(z)}{(z+it)(z-it)}dz.
\end{eqnarray}

We have
\begin{eqnarray}\label{July3a}
&&\frac{1}{\pi i }\int_{c-i\infty}^{c+i\infty}\frac{z{\cal F}_{\tau,t}(z)}{(z+it)(z-it)}dz\nonumber\\
&=&\frac{1}{\pi }\int_{-\infty}^{\infty}\frac{c+iu}{(c+iu)^2+t^2}\cdot\frac{\xi(\tau+\frac{1}{2}+c+iu)\xi(-\tau+\frac{1}{2}+c+iu)}{[(\tau+c+iu)^2-\frac{1}{4}][(-\tau+c+iu)^2-\frac{1}{4}]}du\nonumber\\
&=&\frac{1}{\pi }{\rm Re}\Bigg\{\int_{-\infty}^{\infty}\frac{1}{c+i(u-t)}\cdot\frac{\xi(\tau+\frac{1}{2}+c+iu)\xi(-\tau+\frac{1}{2}+c+iu)}{[(\tau+c+iu)^2-\frac{1}{4}][(-\tau+c+iu)^2-\frac{1}{4}]}du\Bigg\}\nonumber\nonumber\\
&=&\frac{1}{\pi }{\rm Re}\Bigg\{\int_{-\infty}^{\infty}\frac{1}{c+iu}\cdot\frac{\xi(\tau+\frac{1}{2}+c+i(u+t))\xi(-\tau+\frac{1}{2}+c+i(u+t))}{[(\tau+c+i(u+t))^2-\frac{1}{4}][(-\tau+c+i(u+t))^2-\frac{1}{4}]}du\Bigg\}\nonumber\\
&=&\frac{1}{2\pi }\int_{-\infty}^{\infty}\frac{1}{c+iu}\cdot\frac{\xi(\tau+\frac{1}{2}+c+i(u+t))\xi(-\tau+\frac{1}{2}+c+i(u+t))}{[(\tau+c+i(u+t))^2-\frac{1}{4}][(-\tau+c+i(u+t))^2-\frac{1}{4}]}du\nonumber\\
&&+\frac{1}{2\pi }\int_{-\infty}^{\infty}\frac{1}{c+iu}\cdot\frac{\xi(\tau+\frac{1}{2}+c+i(u-t))\xi(-\tau+\frac{1}{2}+c+i(u-t))}{[(\tau+c+i(u-t))^2-\frac{1}{4}][(-\tau+c+i(u-t))^2-\frac{1}{4}]}du\nonumber\\
&=&\frac{1}{2\pi }\int_{-\infty}^{\infty}\left[\int_{1}^{\infty}x^{-(c+iu+1)}dx\right]\frac{\xi(\tau+\frac{1}{2}+c+i(u+t))\xi(-\tau+\frac{1}{2}+c+i(u+t))}{[(\tau+c+i(u+t))^2-\frac{1}{4}][(-\tau+c+i(u+t))^2-\frac{1}{4}]}du\nonumber\\
&&+\frac{1}{2\pi }\int_{-\infty}^{\infty}\left[\int_{1}^{\infty}x^{-(c+iu+1)}dx\right]\frac{\xi(\tau+\frac{1}{2}+c+i(u-t))\xi(-\tau+\frac{1}{2}+c+i(u-t))}{[(\tau+c+i(u-t))^2-\frac{1}{4}][(-\tau+c+i(u-t))^2-\frac{1}{4}]}du\nonumber\\
&=&\frac{1}{2\pi }\int_{-\infty}^{\infty}\left[\int_{1}^{\infty}x^{-[c+i(u-t)+1]}dx\right]\frac{\xi(\tau+\frac{1}{2}+c+iu)\xi(-\tau+\frac{1}{2}+c+iu)}{[(\tau+c+iu)^2-\frac{1}{4}][(-\tau+c+iu)^2-\frac{1}{4}]}du\nonumber\\
&&+\frac{1}{2\pi }\int_{-\infty}^{\infty}\left[\int_{1}^{\infty}x^{-[c+i(u+t)+1]}dx\right]\frac{\xi(\tau+\frac{1}{2}+c+iu)\xi(-\tau+\frac{1}{2}+c+iu)}{[(\tau+c+iu)^2-\frac{1}{4}][(-\tau+c+iu)^2-\frac{1}{4}]}du\nonumber\\
&=&\frac{1}{\pi }\int_{-\infty}^{\infty}\left[\int_{1}^{\infty}x^{-(c+iu+1)}\cos(t\ln x)dx\right]\frac{\xi(\tau+\frac{1}{2}+c+iu)\xi(-\tau+\frac{1}{2}+c+iu)}{[(\tau+c+iu)^2-\frac{1}{4}][(-\tau+c+iu)^2-\frac{1}{4}]}du\nonumber\\
&=&\frac{1}{\pi i}\int_{c-i\infty}^{c+i\infty}\int_{1}^{\infty}x^{-(s+1)}\cos(t\ln x)\cdot\frac{\xi(-\tau+\frac{1}{2}+s)\xi(\tau+\frac{1}{2}+s)}{[(-\tau+s)^2-\frac{1}{4}][(\tau+s)^2-\frac{1}{4}]}dxds\nonumber\\
&=&\int_{1}^{\infty}x^{-1}\cos(t\ln x)\left[\frac{1}{\pi i}\int_{c-i\infty}^{c+i\infty}x^{-s}\cdot\frac{\xi(-\tau+\frac{1}{2}+s)\xi(\tau+\frac{1}{2}+s)}{[(-\tau+s)^2-\frac{1}{4}][(\tau+s)^2-\frac{1}{4}]}ds\right]dx.
\end{eqnarray}
For a function $f$ on $(0,\infty)$, let ${\cal M}f$ denote its Mellin transform. Note that for ${\rm Re}(s)>\frac{1}{2}$,
$$
\{{\cal M}\psi\}(s)=\sum_{n=1}^{\infty}\frac{\Gamma(s)}{(\pi n^2)^s}=\frac{\Gamma(s)\zeta(2s)}{\pi ^{s}}=\frac{\xi(2s)}{s(2s-1)}.
$$
Then, by (\ref{July3a}), we get
\begin{eqnarray}\label{July3b}
&&\frac{1}{\pi i }\int_{c-i\infty}^{c+i\infty}\frac{z{\cal F}_{\tau,t}(z)}{(z+it)(z-it)}dz\nonumber\\
&=&\int_{1}^{\infty}x^{-1}\cos(t\ln x)\left[\frac{1}{4\pi i}\int_{c-i\infty}^{c+i\infty}x^{-s}\{{\cal M}\psi\}\left(\frac{-2\tau+1+2s}{4}\right)\{{\cal M}\psi\}\left(\frac{2\tau+1+2s}{4}\right)ds\right]dx.\nonumber\\
&&
\end{eqnarray}

For $x>0$, define
$$
{\cal L}(x):=\int_{0}^{\infty}x^{\frac{1}{2}}\psi\left(xy\right)\psi\left(xy^{-1}\right)y^{\tau-1}dy.
$$
We have
\begin{eqnarray*}
{\cal L}(x)&=&\int_{0}^{\infty}x^{-\tau+\frac{1}{2}}\psi\left(z\right)\psi\left(x^2z^{-1}\right)z^{\tau-1}dz\\
&=&(x^2)^{\frac{-2\tau+1}{4}}\int_{0}^{\infty}y^{\tau-1}\psi\left(x^2y^{-1}\right)\psi\left(y\right)dy\\
&:=&{\cal G}(x^2).
\end{eqnarray*}
For $x\ge 1$ and any $k\in\mathbb{N}$, we have
\begin{eqnarray*}
{\cal G}(x)&=&x^{\frac{-2\tau+1}{4}}\int_{0}^{\infty}y^{\tau-1}\psi\left(xy^{-1}\right)\psi\left(y\right)dy\\
&=&x^{\frac{-2\tau+1}{4}}\int_{0}^{1}y^{\tau-1}\psi\left(xy^{-1}\right)\psi\left(y\right)dy+x^{\frac{-2\tau+1}{4}}\int_{1}^{\infty}y^{\tau-1}\psi\left(xy^{-1}\right)\psi\left(y\right)dy\\
&=&x^{\frac{-2\tau+1}{4}}\int_{0}^{1}y^{\tau-1}\psi\left(xy^{-1}\right)\psi\left(y\right)dy+x^{\frac{2\tau+1}{4}}\int_{\frac{1}{x}}^{\infty}z^{\tau-1}\psi\left(z^{-1}\right)\psi\left(xz\right)dz\\
&<&x^{\frac{-2\tau+1}{4}}\int_{0}^{1}y^{\tau-1}\psi\left(xy^{-1}\right)\psi\left(y\right)dy+x^{\frac{2\tau+1}{4}}\int_{0}^{\infty}z^{\tau-1}\psi\left(z^{-1}\right)\sum_{n=1}^{\infty}\frac{k!}{(\pi n^2xz)^k}dz.
\end{eqnarray*}
Since $k$ is arbitrary, ${\cal G}(x)$ decays rapidly as $x\rightarrow\infty$. For $0<x<1$, we have
\begin{eqnarray*}
{\cal G}(x)&=&x^{\frac{-2\tau+1}{4}}\int_{0}^{\infty}y^{\tau-1}\psi\left(xy^{-1}\right)\psi\left(y\right)dy\\
&<&x^{\frac{-2\tau+1}{4}}\int_{0}^{\infty}y^{\tau-1}\left(\sum_{n=1}^{\infty}\frac{1}{\pi n^2xy^{-1}}\right)\psi\left(y\right)dy\\
&=&\frac{\pi\cdot x^{\frac{-2\tau-3}{4}}}{6}\int_{0}^{\infty}y^{\tau}\psi\left(y\right)dy.
\end{eqnarray*}
Hence,  by the Mellin inversion theorem (cf. \cite{WM}), we get
\begin{eqnarray}\label{July3c}
{\cal L}(x)&=&\frac{1}{2\pi i}\int_{c-i\infty}^{c+i\infty}x^{-s}\{{\cal M}{\cal L}\}(s)ds\nonumber\\
&=&\frac{1}{4\pi i}\int_{c-i\infty}^{c+i\infty}x^{-s}\{{\cal M}{\cal G}\}\left(\frac{s}{2}\right)ds\nonumber\\
&=&\frac{1}{4\pi i}\int_{c-i\infty}^{c+i\infty}x^{-s}\{{\cal M}\psi\}\left(\frac{-2\tau+1+2s}{4}\right)\{{\cal M}\psi\}\left(\frac{2\tau+1+2s}{4}\right)ds.
\end{eqnarray}
Thus, by (\ref{July3b}) and (\ref{July3c}), we obtain that
\begin{eqnarray*}
\frac{1}{\pi i }\int_{c-i\infty}^{c+i\infty}\frac{z{\cal F}_{\tau,t}(z)}{(z+it)(z-it)}dz=\int_{1}^{\infty}x^{-1}\cos(t\ln x){\cal L}(x)dx,
\end{eqnarray*}
which together with (\ref{July3d}) implies that
\begin{eqnarray}\label{July3e}
&&{\cal F}_{\tau,t}(it)+\frac{[\xi(-2\tau)-\xi(2\tau)]t^2+[(\tau-\frac{1}{2})^2\xi(-2\tau)-(\tau+\frac{1}{2})^2\xi(2\tau)]}{4\tau[t^2+(\tau+\frac{1}{2})^2][t^2+(\tau-\frac{1}{2})^2]}\nonumber\\
&=&\int_{1}^{\infty}x^{-1}\cos(t\ln x){\cal L}(x)dx\nonumber\\
&=&2\int_1^{\infty}\cos(2t\ln x)\int_0^{\infty}y^{\tau-1}\psi(x^2y)\psi(x^2y^{-1})dydx\nonumber\\
&=&\frac{1}{2}\int_{0}^{\infty}\int_{0}^{\infty}\cos\left(\frac{t}{2}\ln (uv)\right)u^{\frac{2\tau-3}{4}}v^{\frac{-2\tau-3}{4}}1_{\{uv>1\}}\psi(u)\psi(v)dudv.
\end{eqnarray}

By (\ref{Rie}) and (\ref{July3e}), we find that for any $\tau>0$ there exists a constant $C_{\tau}$ such that
\begin{eqnarray}\label{Apr4a}
&&4\left|\xi\left(\frac{1}{2}+\tau-it\right)\right|^2\nonumber\\
&=&2\left[\left(t^2+\tau^2+\frac{1}{4}\right)^2-\tau^2\right]\int_{0}^{\infty}\int_{0}^{\infty}\cos\left(\frac{t}{2}\ln (uv)\right)u^{\frac{2\tau-3}{4}}v^{\frac{-2\tau-3}{4}}1_{\{uv>1\}}\psi(u)\psi(v)dudv\nonumber\\
&&-t^2\int_1^{\infty}\left[(1+2\tau)(u^{\tau-\frac{1}{2}}+u^{-\tau-1})+(1-2\tau)(u^{-\tau-\frac{1}{2}}+u^{\tau-1})\right]\psi(u)du\nonumber\\
&&+C_{\tau},\ \ \ \ \forall t\ge0.
\end{eqnarray}
By (\ref{Rie}), we have
\begin{eqnarray}\label{Apr4b}
4\left|\xi\left(\frac{1}{2}+\tau\right)\right|^2=
\left[1+\left(\tau^2-\frac{1}{4}\right)\int_1^{\infty}\left(y^{\frac{2\tau-3}{4}}+y^{\frac{-2\tau-3}{4}}\right)\psi(y)dy\right]^2.
\end{eqnarray}
Comparing (\ref{Apr4a}) with (\ref{Apr4b}), we obtain that
\begin{eqnarray*}
C_{\tau}
&=&\left[1+\left(\tau^2-\frac{1}{4}\right)\int_1^{\infty}\left(y^{\frac{2\tau-3}{4}}+y^{\frac{-2\tau-3}{4}}\right)\psi(y)dy\right]^2\\
&&-2\left(\tau^2-\frac{1}{4}\right)^2\int_{0}^{\infty}\int_{0}^{\infty}u^{\frac{2\tau-3}{4}}v^{\frac{-2\tau-3}{4}}1_{\{uv>1\}}\psi(u)\psi(v)dudv.
\end{eqnarray*}
Therefore, the proof is complete by (\ref{Apr4a}). \hfill \fbox

\begin{rem}
We refer the interested reader to \cite{WS1, WS2} for another proof of the modulus representation (\ref{Apr4v}).  This original proof is more tedious but it indicates how  (\ref{Apr4v}) was discovered.
\end{rem}

\section{Proof of Theorem \ref{thm2}}\setcounter{equation}{0}

We first prove a lemma. For a real number $x$, let $[x]$ denote its integer part.
\begin{lem}\label{July3ab}
For any $\tau\in\mathbb{R}$ and $k\in\{0\}\cup\mathbb{N}$,
\begin{eqnarray*}
&&\int_{0}^{\infty}\int_{0}^{\infty}\left[\frac{1}{2}\ln (uv)\right]^{2k}u^{\frac{2\tau-3}{4}}v^{\frac{-2\tau-3}{4}}1_{\{uv>1\}}\psi(u)\psi(v)dudv\\
&<&\frac{\pi^4[\Gamma(|\tau|+\frac{1}{2})+\Gamma(\frac{1}{2})]([|\tau|]+3)!}{36e}\cdot 2^{2k-\frac{3}{2}}k!.
\end{eqnarray*}
\end{lem}

\noindent {\bf Proof.}\ \ For  $\tau\in\mathbb{R}$, define
$$
{\eta}_{\tau}(y)=\frac{(y+\sqrt{y^2-1})^{\tau}+(y+\sqrt{y^2-1})^{-\tau}}{\sqrt{y^2-1}},\ \ \ \ y>1.
$$
Then, for $k\in\{0\}\cup\mathbb{N}$, we have 
\begin{eqnarray*}
&&\int_{0}^{\infty}\int_{0}^{\infty}\left[\frac{1}{2}\ln (uv)\right]^{2k}u^{\frac{2\tau-3}{4}}v^{\frac{-2\tau-3}{4}}1_{\{uv>1\}}\psi(u)\psi(v)dudv\\
&=&4\int_1^{\infty}(2\ln x)^{2k}\int_0^{\infty}y^{\tau-1}\psi(x^2y)\psi(x^2y^{-1})dydx\\
&=&4\int_1^{\infty}(2\ln x)^{2k}\sum_{m,n=1}^{\infty}\int_{0}^{\infty}u^{\tau-1}e^{-\pi(m^2x^2u+n^2x^2u^{-1})}dudx\\
&=&4\int_1^{\infty}(2\ln x)^{2k}\left\{\sum_{m,n=1}^{\infty}\int_{\frac{n}{m}}^{\infty}u^{\tau-1}e^{-\pi(m^2x^2u+n^2x^2u^{-1})}du+\sum_{m,n=1}^{\infty}\int_0^{\frac{n}{m}}u^{\tau-1}e^{-\pi(n^2x^2u^{-1}+m^2x^2u)}du\right\}dx\\
&=&4\int_1^{\infty}(2\ln x)^{2k}\left\{\sum_{m,n=1}^{\infty}\int_{\frac{n}{m}}^{\infty}u^{\tau-1}e^{-\pi(m^2x^2u+n^2x^2u^{-1})}du+\sum_{m,n=1}^{\infty}\int_{\frac{m}{n}}^{\infty}v^{-\tau-1}e^{-\pi(n^2x^2v+m^2x^2v^{-1})}dv\right\}dx\\
&=&4\int_1^{\infty}(2\ln x)^{2k}\left\{\sum_{m,n=1}^{\infty}\left(\frac{n}{m}\right)^{\tau}\int_1^{\infty}z^{\tau-1}e^{-\pi mnx^2(z+z^{-1})}dz+\sum_{m,n=1}^{\infty}\left(\frac{n}{m}\right)^{\tau}\int_1^{\infty}z^{-\tau-1}e^{-\pi mnx^2(z+z^{-1})}dz\right\}dx\\
&=&4\int_1^{\infty}\int_1^{\infty}(2\ln x)^{2k}\sum_{m,n=1}^{\infty}\left(\frac{n}{m}\right)^{\tau}\int_1^{\infty} e^{-2\pi mnx^2y}{\eta}_{\tau}(y)dydx\\
&=&2\sum_{m,n=1}^{\infty}\left(\frac{n}{m}\right)^{\tau}\int_1^{\infty} \int_1^{\infty}\frac{x^{-\frac{1}{2}}(\ln x)^{2k}e^{-2\pi mn xy}[(y+\sqrt{y^2-1})^{\tau}+(y+\sqrt{y^2-1})^{-\tau}]}{\sqrt{y^2-1}}dxdy\\
&<&4\Bigg\{\sum_{m\le n}^{\infty}\left(\frac{n}{m}\right)^{|\tau|}e^{- 4mn }\Bigg\}\int_1^{\infty} \frac{[(2y)^{|\tau|}+1]e^{- y}}{\sqrt{y^2-1}}dy\int_1^{\infty}x^{-\frac{1}{2}}(2x^{\frac{1}{2}})^{2k}e^{- x}dx\\
&<&4\Bigg\{\sum_{m,n=1}^{\infty}\left(\frac{n}{m}\right)^{|\tau|}\cdot\frac{([|\tau|]+3)!}{(4mn)^{[|\tau|]+3}}\Bigg\}\cdot\frac{1}{\sqrt{2}e}\int_0^{\infty} \frac{2^{2|\tau|+1}(y^{|\tau|}+1)e^{- y}}{\sqrt{y}}dy\cdot 2^{2k}\int_0^{\infty}x^{k}e^{-x}dx\\
&<&\frac{\pi^4[\Gamma(|\tau|+\frac{1}{2})+\Gamma(\frac{1}{2})]([|\tau|]+3)!}{36e}\cdot 2^{2k-\frac{3}{2}}k!.
\end{eqnarray*}
The proof is complete.\hfill \fbox

\vskip 0.3cm
\noindent {\bf Proof of Theorem \ref{thm2}.}\ \ To simplify notation, define
\begin{eqnarray}\label{Nov2bv}
F(\tau,t):=4\left|\xi\left(\frac{1}{2}+\tau-it\right)\right|^2,\ \ \ \ \tau,t\in\mathbb{R}.
\end{eqnarray} 
By (\ref{Apr4v}) and Lemma \ref{July3ab}, we  get
\begin{eqnarray}\label{July3ae}
&&\partial_{\tau}F(\tau,t)\nonumber\\
&=&8\tau\left[t^2-\left(\frac{1}{4}-\tau^2\right)\right]\int_{0}^{\infty}\int_{0}^{\infty}\cos\left(\frac{t}{2}\ln (uv)\right)u^{\frac{2\tau-3}{4}}v^{\frac{-2\tau-3}{4}}1_{\{uv>1\}}\psi(u)\psi(v)dudv\nonumber\\
&&+\left[\left(t^2+\frac{1}{4}+\tau^2\right)^2-\tau^2\right]\int_{0}^{\infty}\int_{0}^{\infty}\cos\left(\frac{t}{2}\ln (uv)\right)\ln\left(\frac{u}{v}\right)u^{\frac{2\tau-3}{4}}v^{\frac{-2\tau-3}{4}}1_{\{uv>1\}}\psi(u)\psi(v)dudv\nonumber\\
&&-2t^2\int_1^{\infty}\left[(u^{\tau-\frac{1}{2}}+u^{-\tau-1})-(u^{-\tau-\frac{1}{2}}+u^{\tau-1})\right]\psi(u)du\nonumber\\
&&-t^2\int_1^{\infty}\left[(1+2\tau)(u^{\tau-\frac{1}{2}}-u^{-\tau-1})+(1-2\tau)(u^{\tau-1}-u^{-\tau-\frac{1}{2}})\right](\ln u)\psi(u)du\nonumber\\
&&+8\tau\left(\frac{1}{4}-\tau^2\right)\int_{0}^{\infty}\int_{0}^{\infty}u^{\frac{2\tau-3}{4}}v^{\frac{-2\tau-3}{4}}1_{\{uv>1\}}\psi(u)\psi(v)dudv\nonumber\\
&&-\left(\frac{1}{4}-\tau^2\right)^2\int_{0}^{\infty}\int_{0}^{\infty}\ln\left(\frac{u}{v}\right)u^{\frac{2\tau-3}{4}}v^{\frac{-2\tau-3}{4}}1_{\{uv>1\}}\psi(u)\psi(v)dudv\nonumber\\
&&+\left[1-\left(\frac{1}{4}-\tau^2\right)\int_1^{\infty}\left(u^{\frac{2\tau-3}{4}}+u^{\frac{-2\tau-3}{4}}\right)\psi(u)du\right]\nonumber\\
&&\quad\cdot\left[4\tau\int_1^{\infty}\left(u^{\frac{2\tau-3}{4}}+u^{\frac{-2\tau-3}{4}}\right)\psi(u)du-\left(\frac{1}{4}-\tau^2\right)\int_1^{\infty}\left(u^{\frac{2\tau-3}{4}}-u^{\frac{-2\tau-3}{4}}\right)(\ln u)\psi(u)du\right].\nonumber\\
&&
\end{eqnarray}

Note that
$$
\cos x=\sum_{k=0}^{\infty}\frac{(-1)^k}{(2k)!}x^{2k},\ \ \ \ x\in\mathbb{R},
$$
and
$$
\sum_{k=0}^{2n-1}\frac{(-1)^k}{(2k)!}x^{2k}\le\cos x\le \sum_{k=0}^{2n-2}\frac{(-1)^k}{(2k)!}x^{2k},\ \ \ \ x\in\mathbb{R},\ n\in\mathbb{N}.
$$
Then, for any $n\in\mathbb{N}$, 
\begin{eqnarray}\label{Nov2a}
&&\sum_{k=0}^{2n-1}\frac{(-1)^kt^{2k}}{(2k)!\cdot 2^{2k}}\int_{0}^{\infty}\int_{0}^{\infty}\left[\ln (uv)\right]^{2k}u^{\frac{2\tau-3}{4}}v^{\frac{-2\tau-3}{4}}1_{\{uv>1\}}\psi(u)\psi(v)dudv\nonumber\\
&<&\int_{0}^{\infty}\int_{0}^{\infty}\cos\left(\frac{t}{2}\ln (uv)\right)u^{\frac{2\tau-3}{4}}v^{\frac{-2\tau-3}{4}}1_{\{uv>1\}}\psi(u)\psi(v)dudv\nonumber\\
&<&\sum_{k=0}^{2n-2}\frac{(-1)^kt^{2k}}{(2k)!\cdot 2^{2k}}\int_{0}^{\infty}\int_{0}^{\infty}\left[\ln (uv)\right]^{2k}u^{\frac{2\tau-3}{4}}v^{\frac{-2\tau-3}{4}}1_{\{uv>1\}}\psi(u)\psi(v)dudv,
\end{eqnarray}
and
\begin{eqnarray}\label{Nov2b}
&&\int_{0}^{\infty}\int_{0}^{\infty}\cos\left(\frac{t}{2}\ln (uv)\right)\ln\left(\frac{u}{v}\right)u^{\frac{2\tau-3}{4}}v^{\frac{-2\tau-3}{4}}1_{\{uv>1\}}\psi(u)\psi(v)dudv\nonumber\\
&=&\int_{0}^{\infty}\int_{0}^{\infty}\cos\left(\frac{t}{2}\ln (uv)\right)\ln\left(\frac{u}{v}\right)\left(\frac{u}{v}\right)^{\frac{\tau}{2}}(uv)^{-\frac{3}{4}}1_{\{uv>1,u>v\}}\psi(u)\psi(v)dudv\nonumber\\
&&+\int_{0}^{\infty}\int_{0}^{\infty}\cos\left(\frac{t}{2}\ln (uv)\right)\ln\left(\frac{u}{v}\right)\left(\frac{u}{v}\right)^{\frac{\tau}{2}}(uv)^{-\frac{3}{4}}1_{\{uv>1,u<v\}}\psi(u)\psi(v)dudv\nonumber\\
&=&\int_{0}^{\infty}\int_{0}^{\infty}\cos\left(\frac{t}{2}\ln (uv)\right)\ln\left(\frac{u}{v}\right)\left[\left(\frac{u}{v}\right)^{\frac{\tau}{2}}-\left(\frac{v}{u}\right)^{\frac{\tau}{2}}\right](uv)^{-\frac{3}{4}}1_{\{uv>1,u>v\}}\psi(u)\psi(v)dudv\nonumber\\
&<&\sum_{k=0}^{2n-2}\frac{(-1)^kt^{2k}}{(2k)!\cdot 2^{2k}}\int_{0}^{\infty}\int_{0}^{\infty}\left[\ln (uv)\right]^{2k}\ln\left(\frac{u}{v}\right)\left[\left(\frac{u}{v}\right)^{\frac{\tau}{2}}-\left(\frac{v}{u}\right)^{\frac{\tau}{2}}\right](uv)^{-\frac{3}{4}}1_{\{uv>1,u>v\}}\psi(u)\psi(v)dudv\nonumber\\
&=&\sum_{k=0}^{2n-2}\frac{(-1)^kt^{2k}}{(2k)!\cdot 2^{2k}}\int_{0}^{\infty}\int_{0}^{\infty}\left[\ln (uv)\right]^{2k}\ln\left(\frac{u}{v}\right)u^{\frac{2\tau-3}{4}}v^{\frac{-2\tau-3}{4}}1_{\{uv>1\}}\psi(u)\psi(v)dudv.
\end{eqnarray}
Thus, by (\ref{July3ae})--(\ref{Nov2b}), we deduce that for any $t\in\mathbb{R}$, $\tau\in(0,\infty)$ and $n\in\mathbb{N}$,
\begin{eqnarray}\label{Nov2s}
&&\partial_{\tau}F(\tau,t)\nonumber\\
&<&8\tau t^2\sum_{k=0}^{2n-2}\frac{(-1)^kt^{2k}}{(2k)!\cdot 2^{2k}}\int_{0}^{\infty}\int_{0}^{\infty}\left[\ln (uv)\right]^{2k}u^{\frac{2\tau-3}{4}}v^{\frac{-2\tau-3}{4}}1_{\{uv>1\}}\psi(u)\psi(v)dudv\nonumber\\
&&-8\tau\left(\frac{1}{4}-\tau^2\right)\sum_{k=0}^{2n-1}\frac{(-1)^kt^{2k}}{(2k)!\cdot 2^{2k}}\int_{0}^{\infty}\int_{0}^{\infty}\left[\ln (uv)\right]^{2k}u^{\frac{2\tau-3}{4}}v^{\frac{-2\tau-3}{4}}1_{\{uv>1\}}\psi(u)\psi(v)dudv\nonumber\\
&&+\left[t^4+2\left(\frac{1}{4}+\tau^2\right)t^2+\left(\frac{1}{4}-\tau^2\right)^2\right]\nonumber\\
&&\ \ \ \ \cdot\sum_{k=0}^{2n-2}\frac{(-1)^kt^{2k}}{(2k)!\cdot 2^{2k}}\int_{0}^{\infty}\int_{0}^{\infty}\left[\ln (uv)\right]^{2k}\ln\left(\frac{u}{v}\right)u^{\frac{2\tau-3}{4}}v^{\frac{-2\tau-3}{4}}1_{\{uv>1\}}\psi(u)\psi(v)dudv\nonumber\\
&&-2t^2\int_1^{\infty}\left[(u^{\tau-\frac{1}{2}}+u^{-\tau-1})-(u^{-\tau-\frac{1}{2}}+u^{\tau-1})\right]\psi(u)du\nonumber\\
&&-t^2\int_1^{\infty}\left[(1+2\tau)(u^{\tau-\frac{1}{2}}-u^{-\tau-1})+(1-2\tau)(u^{\tau-1}-u^{-\tau-\frac{1}{2}})\right](\ln u)\psi(u)du\nonumber\\
&&+8\tau\left(\frac{1}{4}-\tau^2\right)\int_{0}^{\infty}\int_{0}^{\infty}u^{\frac{2\tau-3}{4}}v^{\frac{-2\tau-3}{4}}1_{\{uv>1\}}\psi(u)\psi(v)dudv\nonumber\\
&&-\left(\frac{1}{4}-\tau^2\right)^2\int_{0}^{\infty}\int_{0}^{\infty}\ln\left(\frac{u}{v}\right)u^{\frac{2\tau-3}{4}}v^{\frac{-2\tau-3}{4}}1_{\{uv>1\}}\psi(u)\psi(v)dudv\nonumber\\
&&+\left[1-\left(\frac{1}{4}-\tau^2\right)\int_1^{\infty}\left(u^{\frac{2\tau-3}{4}}+u^{\frac{-2\tau-3}{4}}\right)\psi(u)du\right]\nonumber\\
&&\quad\cdot\left[4\tau\int_1^{\infty}\left(u^{\frac{2\tau-3}{4}}+u^{\frac{-2\tau-3}{4}}\right)\psi(u)du-\left(\frac{1}{4}-\tau^2\right)\int_1^{\infty}\left(u^{\frac{2\tau-3}{4}}-u^{\frac{-2\tau-3}{4}}\right)(\ln u)\psi(u)du\right]\nonumber\\
&=&f_{\tau,n}(t).
\end{eqnarray}
According to \cite{SD}, the RH is equivalent to the statement that $\partial_{\tau}F(\tau,t)\ge 0$ for all $t\in\mathbb{R}$ and $\tau\in(0,\frac{1}{2})$. Therefore, the proof is complete by (\ref{Nov2s}), Lemma \ref{July3ab} and approximation.\hfill \fbox

\begin{rem}\label{remNov} Note that for any $t\in\mathbb{R}$, $\partial_{\tau}F(0,t)=0$ and $\partial_{\tau}F(\tau,t)\ge0$ for $\tau\ge \frac{1}{2}$ (cf. \cite{SD}). Then, $f_{\tau,n}(t)>0$ for any $t\in\mathbb{R}$, $\tau\in\{0\}\cup[\frac{1}{2},\infty)$ and $n\in\mathbb{N}$.
\end{rem}

\section{Proofs of Theorem \ref{pro99} and Corollary \ref{cor111}}\setcounter{equation}{0}

We first prove a lemma.
\begin{lem}\label{lem41}
Let $\tau\in(0,\frac{1}{2}]$, $n\in\mathbb{N}$ and $k\in\{0\}\cup\mathbb{N}$. Then,

\noindent (i) $a_{\tau,n}(2n)>0$ and $a_{\tau,n}(2n)\uparrow a_{\frac{1}{2},n}(2n)$ as $\tau\uparrow \frac{1}{2}$. 

\noindent (ii) $\sup_{\tau\in(0,\frac{1}{2}]}|a_{\tau,n}(2n-1)|<\infty$ and $\sup_{\tau\in(0,\frac{1}{2}]}|a_{\tau}(k)|<\infty$. 

\noindent (iii)  $a_{\tau}(0)=\partial_{\tau}F(\tau,0)>0$, where the function $F$ is given by (\ref{Nov2bv}).

\end{lem}

\noindent {\bf Proof.}\ \ (i) We have
\begin{eqnarray*}
a_{\tau,n}(2n)&=&\frac{1}{(4n-4)!\cdot 2^{4n-4}}\int_{0}^{\infty}\int_{0}^{\infty}\left[\ln (uv)\right]^{4n-4}\ln\left(\frac{u}{v}\right)\left(\frac{u}{v}\right)^{\frac{\tau}{2}}(uv)^{-\frac{3}{4}}1_{\{uv>1,u>v\}}\psi(u)\psi(v)dudv\\
&&+\frac{1}{(4n-4)!\cdot 2^{4n-4}}\int_{0}^{\infty}\int_{0}^{\infty}\left[\ln (uv)\right]^{4n-4}\left(\frac{u}{v}\right)^{\frac{\tau}{2}}(uv)^{-\frac{3}{4}}1_{\{uv>1,u<v\}}\psi(u)\psi(v)dudv\\
&=&\frac{1}{(4n-4)!\cdot 2^{4n-4}}\int_{0}^{\infty}\int_{0}^{\infty}\left[\ln (uv)\right]^{4n-4}\ln\left(\frac{u}{v}\right)\left[\left(\frac{u}{v}\right)^{\frac{\tau}{2}}-\left(\frac{v}{u}\right)^{\frac{\tau}{2}}\right](uv)^{-\frac{3}{4}}1_{\{uv>1,u>v\}}\psi(u)\psi(v)dudv.
\end{eqnarray*}
Then, the assertions hold. 

\noindent (ii) We can prove the two assertions by following the similar argument used in (i).

\noindent (iii) We have
\begin{eqnarray*}
\int_1^{\infty}\left(u^{\frac{2\tau-3}{4}}+u^{\frac{-2\tau-3}{4}}\right)\psi(u)du
&<&2\sum_{n=1}^{\infty}\int_1^{\infty}e^{-\pi n^2 u}du\\
&=&\sum_{n=1}^{\infty}\frac{2}{\pi n^2e^{\pi n^2}}\\
&<&\frac{\pi}{3e^{\pi }}\\
&\approx&0.04525351,
\end{eqnarray*}
and
\begin{eqnarray*}
&&4\tau\int_1^{\infty}\left(u^{\frac{2\tau-3}{4}}+u^{\frac{-2\tau-3}{4}}\right)\psi(u)du-\left(\frac{1}{4}-\tau^2\right)\int_1^{\infty}\left(u^{\frac{2\tau-3}{4}}-u^{\frac{-2\tau-3}{4}}\right)(\ln u)\psi(u)du\\
&=&4\tau\int_1^{\infty}\left(u^{\frac{2\tau-3}{4}}+u^{\frac{-2\tau-3}{4}}\right)\psi(u)du-\left(\frac{1}{4}-\tau^2\right)\int_1^{\infty}\left\{\int_{\frac{-2\tau-3}{4}}^{\frac{2\tau-3}{4}}u^x(\ln u)dx\right\}(\ln u)\psi(u)du\\
&>&4\tau\int_1^{\infty}u^{\frac{2\tau-3}{4}}\psi(u)du-\frac{\tau}{4}\int_1^{\infty}u^{\frac{2\tau-3}{4}}(\ln u)^2\psi(u)du\\
&>&\frac{\tau}{4}\left\{15\int_1^{\infty}u^{\frac{2\tau-3}{4}}\psi(u)du-\int_e^{\infty}u^{\frac{2\tau-3}{4}}(\ln u)^2\psi(u)du\right\}\\
&>&\frac{\tau}{4}\left\{15\int_1^{\infty}u^{-\frac{3}{4}}\psi(u)du-\int_e^{\infty}u^{\frac{3}{2}}\psi(u)du\right\}\\
&=&\frac{\tau}{4}\left\{\sum_{n=1}^{\infty}\left[15\int_1^{\infty}u^{-\frac{3}{4}}e^{-\pi n^2 u}du-\int_e^{\infty}u^{\frac{3}{2}}e^{-\pi n^2 u}du\right]\right\}\\
&>&\frac{\tau}{4}\left\{\sum_{n=1}^{\infty}e^{-2\pi n^2 }\left[15\cdot2^{-\frac{3}{4}}-\int_{e-2}^{\infty}(v+2)^{\frac{3}{2}}e^{-\pi n^2 v}dv\right]\right\}\\
&>&\frac{\tau}{4}\left\{\sum_{n=1}^{\infty}e^{-2\pi n^2 }\left[15\cdot2^{-\frac{3}{4}}-\int_{0}^{\infty}(v+2)^{2}e^{-\pi v}dv\right]\right\}\\
&=&\frac{\tau}{4}\left\{\sum_{n=1}^{\infty}e^{-2\pi n^2 }\left[15\cdot2^{-\frac{3}{4}}-\frac{1}{\pi}\left(\frac{2}{\pi^2}+\frac{4}{\pi}+4\right)\right]\right\}\\
&\approx&\frac{\tau}{4}\left\{7.176026\sum_{n=1}^{\infty}e^{-2\pi n^2 }\right\}.
\end{eqnarray*}
Then,  $\partial_{\tau}F(\tau,0)=a_{\tau}(0)>0$ by (\ref{July4b}) and (\ref{July3ae}).\hfill \fbox

\begin{rem}
We would like to point out that the sign of each term appearing in the definitions of $a_{\tau}(k)$ and $a_{\tau,n}(2n-1)$, $k,n\in\mathbb{N}$, $\tau\in(0,\frac{1}{2})$, is either positive or negative. For example, the last two terms appearing in the definition of $a_{\tau}(1)$ are negative for $\tau\in(0,\frac{1}{2})$. In fact, for $u>1$, we have
\begin{eqnarray*}
(u^{\tau-\frac{1}{2}}+u^{-\tau-1})-(u^{-\tau-\frac{1}{2}}+u^{\tau-1})
&=&u^{-\tau-1}(u^{2\tau}-1)(u^{\frac{1}{2}}-1)\\
&>&0,
\end{eqnarray*}
for $\tau\in[\frac{1}{4},\frac{1}{2})$,
\begin{eqnarray*}
&&(1+2\tau)(u^{\tau-\frac{1}{2}}-u^{-\tau-1})+(1-2\tau)(u^{\tau-1}-u^{-\tau-\frac{1}{2}})\\
&\ge&(1+2\tau)(u^{\tau-\frac{1}{2}}-u^{-\tau-1})\\
&>&0,
\end{eqnarray*}
and for $\tau\in(0,\frac{1}{4})$,
\begin{eqnarray*}
&&(1+2\tau)(u^{\tau-\frac{1}{2}}-u^{-\tau-1})+(1-2\tau)(u^{\tau-1}-u^{-\tau-\frac{1}{2}})\\
&=&(1+2\tau)(u^{\tau-\frac{1}{2}}-u^{-\tau-1})-(1-2\tau)(u^{-\tau-\frac{1}{2}}-u^{\tau-1})\\
&>&(1+2\tau)(u^{-\tau-\frac{1}{2}}-u^{\tau-1})-(1-2\tau)(u^{-\tau-\frac{1}{2}}-u^{\tau-1})\\
&=&4\tau(u^{-\tau-\frac{1}{2}}-u^{\tau-1})\\
&>&0.
\end{eqnarray*}
\end{rem}

\noindent {\bf Proof of Theorem \ref{pro99}.}\ \ By Theorem \ref{thm2}, we can prove the necessity. Below we prove the sufficiency.  Suppose that $n_k\uparrow\infty$ and (\ref{Nov1a}) holds for any $k\in\mathbb{N}$ and $\tau\in(0,\frac{1}{2})$. Fix a $k_0\in\mathbb{N}$. We will show that $f_{\tau,n_{k_0}}(t)>0$ for any $t\in\mathbb{R}$ and $\tau\in(0,\frac{1}{2})$. Once this is proved, the sufficiency is obtained by Theorem \ref{thm2}.

\noindent (a) First, we consider the case that ${\rm Discr}(f_{\tau,n_{k_0}})\not=0$ for any $\tau\in(0,\frac{1}{2})$, i.e., the polynomial $f_{\tau,n_{k_0}}$ has no multiple roots for any $\tau\in(0,\frac{1}{2})$. Note that the Hermite form $H_1(f_{\tau,n_{k_0}})$ of $f_{\tau,n_{k_0}}$ is a real-valued symmetric  matrix (cf. \cite{P}). Then, all eigenvalues of $H_1(f_{\tau,n_{k_0}})$ are real numbers.

Assume the contrary that ${\cal N}(f_{\tau_0,n_{k_0}})\ge1$ for some $\tau_0\in(0,\frac{1}{2})$. Then, the signature of  $H_1(f_{\tau_0,n_{k_0}})$ of $f_{\tau_0,n_{k_0}}$ is $\ge 1$. Since $f_{\tau_0,n_{k_0}}$ has no multiple roots, the rank of $H_1(f_{\tau_0,n_{k_0}})$ is equal to $4n$ (cf. \cite{P}), which implies that 0 is not an eigenvalue of $H_1(f_{\tau_0,n_{k_0}})$. Hence, by continuity of the roots of a (characteristic) polynomial, there exists an open interval containing $\tau_0$ such that for any point $\tau$ of this interval the signature of Hermite form $H_1(f_{\tau,n_{k_0}})$ is $\ge 1$.

Denote by $I_{n_{k_0}}=(a_{n_{k_0}}, b_{n_{k_0}})$ the largest open interval containing $\tau_0$ such that the signature of $H_1(f_{\tau,n_{k_0}})$ is $\ge 1$ for any $\tau\in I_{n_{k_0}}$. We have $b_{n_{k_0}}\in (0,\frac{1}{2}]$. By continuity of the roots of a polynomial, we find that 0 must be an eigenvalue of $H_1(f_{b_{n_{k_0}},n_{k_0}})$. Since the rank of $H_1(f_{b_{n_{k_0}},n_{k_0}})$ is equal to the number of distinct complex roots of $f_{b_{n_{k_0}},n_{k_0}}$, we deduce that $b_{n_{k_0}}=\frac{1}{2}$. Thus, ${\cal N}(f_{\tau,n_{k_0}})\ge 1$  for any $\tau\in(a_{n_{k_0}}, \frac{1}{2})$. By continuity and Lemma \ref{lem41} (i) and (ii), this implies that ${\cal N}(f_{\frac{1}{2},n_{k_0}})\ge 1$. Hence, $\partial_{\tau}F(\frac{1}{2},t)<0$ for some $t\in\mathbb{R}$ by (\ref{Nov2s}). 
We have arrived at a contradiction by \cite[Theorem 1]{SD}.

\noindent (b) Next, we consider the case that
$$
A_{k_0}:=\left\{\tau\in\left(0,\frac{1}{2}\right):{\rm Discr}(f_{\tau,n_{k_0}})\not=0\right\}\not=\emptyset.
$$
We will show that $f_{\tau,n_{k_0}}(t)>0$ for any $t\in\mathbb{R}$ and $\tau\in A_{k_0}$.

Assume the contrary that  ${\cal N}(f_{\tau_0,n_{k_0}})\ge 1$ for some $\tau_0\in A_{k_0}$. Similar to (a), we can show that there exists an open interval containing $\tau_0$ such that for any point $\tau$ of this interval ${\rm Discr}(f_{\tau,n_{k_0}})\not=0$ and the signature of Hermite form $H_1(f_{\tau,n_{k_0}})$ is $\ge 1$.
Denote by $I_{n_{k_0}}=(a_{n_{k_0}}, b_{n_{k_0}})$ the largest open interval containing $\tau_0$ such that ${\rm Discr}(f_{\tau,n_{k_0}})\not=0$ and the signature of $H_1(f_{\tau,n_{k_0}})$ is $\ge 1$ for any $\tau\in I_{n_{k_0}}$. We have $b_{n_{k_0}}\in (0,\frac{1}{2}]$. By continuity of the roots of a polynomial, we find that 0 must be an eigenvalue of $H_1(f_{b_{n_{k_0}},n_{k_0}})$. Since the rank of $H_1(f_{b_{n_{k_0}},n_{k_0}})$ is equal to the number of distinct complex roots of $f_{b_{n_{k_0}},n_{k_0}}$, $b_{n_{k_0}}\in (0,\frac{1}{2})\cap A^c_{k_0}$ or  $b_{n_{k_0}}=\frac{1}{2}$. Further, by continuity and Lemma \ref{lem41} (i) and (ii), we deduce that ${\cal N}(f_{b_{n_{k_0}},n_{k_0}})\ge 1$ since ${\cal N}(f_{\tau,n_{k_0}})\ge1$ for any $\tau\in(a_{n_{k_0}}, b_{n_{k_0}})$. We have arrived at a contradiction by (\ref{Nov1a}), \cite[Theorem 1]{SD} and (\ref{Nov2s}).\hfill \fbox

\vskip 0.3cm
\noindent {\bf Proof of Corollary \ref{cor111}.}\ \ (i) follows from Theorem \ref{thm2} and a basic property of the discriminant of a polynomial with real coefficients. 

\noindent (ii) By (\ref{Nov4a}), for any $\varepsilon\in(0,\frac{1}{2})$, there exists $k_{\varepsilon}\in\mathbb{N}$ such that 
${\rm Discr}(f_{\tau,n_{k}})\not=0$ for any $\tau\in(\varepsilon,\frac{1}{2})$ and $k\ge k_{\varepsilon}$. Then, following the similar argument used in the proof of Theorem \ref{pro99}, we can show that $f_{\tau,n_{k}}(t)>0$ for any $t\in\mathbb{R}$, $\tau\in(\varepsilon,\frac{1}{2})$  and $k\ge k_{\varepsilon}$. Thus, $\partial_{\tau}F(\tau,t)\ge 0$ for any $t\in\mathbb{R}$ and $\tau\in(\varepsilon,\frac{1}{2})$ by (\ref{Nov2s}). Hence, the RH is true since $\varepsilon\in(0,\frac{1}{2})$ is arbitrary. 

\noindent (iii) We have
\begin{eqnarray}\label{Oct29a}
{\rm Discr}(f_{\tau,1})&=&256[a_{\tau,1}(2)]^3[a_{\tau}(0)]^3-128[a_{\tau,1}(2)]^2[a_{\tau,1}(1)]^2[a_{\tau}(0)]^2+16a_{\tau,1}(2)[a_{\tau,1}(1)]^4a_{\tau}(0)\nonumber\\
&=&16a_{\tau,1}(2)a_{\tau}(0)\left\{[a_{\tau,1}(1)]^2-4a_{\tau,1}(2)a_{\tau}(0)\right\}^2.
\end{eqnarray}
Note that, by Lemma \ref{lem41}, $a_{\tau,1}(2)>0$ and $a_{\tau}(0)>0$ for any $\tau\in(0,\frac{1}{2})$. Then, the assertions hold. 

\noindent (iv) ``$\Longrightarrow$"\ \ Assume that $a_{\tau,1}(1)=-2\sqrt{a_{\tau,1}(2)a_{\tau}(0)}$ for some $\tau\in(0,\frac{1}{2})$. Then, we have
$$
f_{\tau,1}\left(\left(\frac{a_{\tau}(0)}{a_{\tau,1}(2)}\right)^{\frac{1}{4}}\right)=0.
$$
Further, the last assertion holds as a consequence of Theorem \ref{thm2}. 

\noindent ``$\Longleftarrow$"\ \ Assume that $a_{\tau,1}(1)\not=-2\sqrt{a_{\tau,1}(2)a_{\tau}(0)}$ for any $\tau\in(0,\frac{1}{2})$. By (\ref{Oct29a}), we deduce that for any $\tau\in(0,\frac{1}{2})$,
$$
{\rm Discr}(f_{\tau,1})=0\Longrightarrow a_{\tau,1}(1)=2\sqrt{a_{\tau,1}(2)a_{\tau}(0)}>0\Longrightarrow {\cal N}(f_{\tau,1})=0.
$$
Therefore, $f_{\tau,1}(t)>0$ for any $t\in\mathbb{R}$ and $\tau\in(0,\frac{1}{2})$ by the proof of Theorem \ref{pro99}.\hfill \fbox

\begin{rem}
According to Theorem \ref{pro99} and Corollary \ref{cor111}, it is important to study zeros of the functions $\tau\mapsto {\rm Discr}(f_{\tau,n})$, $n\in\mathbb{N}$.
\end{rem}


\begin{thebibliography}{1234}



\bibitem{E}  H.M. Edwards. Riemann's Zeta Function. Academic Press (1974).


\bibitem{P} P.A. Parrilo. Algebraic techniques and semidefinite optimization, Lecture 5. https://ocw.mit.edu/courses/6-972-algebraic-techniques-and-semidefinite-optimization-spring-2006/pages/lecture-notes/ (2010).

\bibitem{R} B. Riemann.  Ueber die Anzahl der Primzahlen unter einer gegebenen Gr\"osse. Monatsberichte der Berliner Akademie (1859).


\bibitem{SD} J. Sondow and C. Dumitrescu. A monotonicity property of Riemann's xi function and a reformulation of the Riemann hypothesis. Period. Math. Hungar. 60,  37-40 (2010).

\bibitem{WS1} W. Sun.  An equivalent inequality for the Riemann hypothesis. arXiv:2312.12760v5 (2024).

\bibitem{WS2} W. Sun.  Modulus representation of the Riemann $\xi$ function. arXiv:2404.06642v1 (2024).

\bibitem{WM}  Wikipedia. Mellin transform. https://en.wikipedia.org/wiki/Mellin\_transform.

\bibitem{W}  Wikipedia. Riemann hypothesis. https://en.wikipedia.org/wiki/Riemann\_hypothesis.

\bibitem{Y1}  L. Yang, X. Hou and  Z.B. Zeng. A complete discrimination system for polynomials. Sci. China, Ser. E 39, 628-646 (1996).

\bibitem{Y2}  L. Yang, J. Zhang and X. Hou. Nonlinear Algebraic Equation System and Automated Theorem Proving. Shanghai Scientific and
Technological Education Publishing House (1996).


\end{thebibliography}
\end{document}